\documentclass[10pt]{article}
\usepackage{amsmath}
\usepackage{amssymb}
\usepackage{amsthm}

\usepackage{tikz}
\usetikzlibrary{matrix,arrows,decorations.pathmorphing}
\usepackage{enumitem}

\newtheorem{definition}{Definition}[section]
\newtheorem{theorem}[definition]{Theorem}
\newtheorem{example}[definition]{Example}
\newtheorem{lemma}[definition]{Lemma}

\newtheorem{remark}[definition]{Remark}

\newtheorem{proposition}[definition]{Proposition}

\newtheorem{notation}[definition]{Notation}
\def\N{{\mathbb N}}
\def\F{{\mathbb F}}
\def\Z{{\mathbb Z}}

\def\H{{\mathbb H}}
\def\I{{\mathbb I}}

\def\Span{{\rm Span}\ }

\def\PSL2{PSL_2 (\mathbb{Z})}
\def\SL2{SL_2 (\mathbb{Z})}

\def\dirsum{\ \ \ \ \ \ \ \ \ \ ({\mbox{direct sum}})}

\def\white{\Delta}
\def\genset{\mathcal{X}}

\def\monoid{\left<\genset\right>}
\def\lbrack{\left[}
\def\rbrack{\right]}
\def\rev{^*}
\def\LieABC{L}
\def\freeLie{\mathcal{L}}
\def\ad{{\rm ad}\ }

\def\Hfour{\lbrack B , A\rbrack}
\def\Hfive{\lbrack C , A\rbrack}
\def\Hsix{\lbrack C , B\rbrack}

\def\Hseven{\lbrack \Hfour , A\rbrack}
\def\Height{\lbrack \Hfive , A\rbrack}
\def\Hnine{\lbrack \Hfour , B\rbrack}
\def\Hten{\lbrack \Hfive , B\rbrack}
\def\Heleven{\lbrack \Hsix , B\rbrack}
\def\Htwelve{\lbrack \Hfour , C\rbrack}
\def\Hthirteen{\lbrack \Hfive , C\rbrack}
\def\Hfourteen{\lbrack \Hsix , C\rbrack}

\def\Hfifteen{\lbrack \Hseven , A\rbrack}
\def\Hsixteen{\lbrack \Height , A\rbrack}
\def\Hseventeen{\lbrack \Hseven , B\rbrack}
\def\Heighteen{\lbrack \Height , B\rbrack}
\def\Hnineteen{\lbrack \Hnine , B\rbrack}
\def\Htwenty{\lbrack \Hten , B\rbrack}
\def\Htwentyone{\lbrack \Heleven , B\rbrack}
\def\Htwentytwo{\lbrack \Hseven , C\rbrack}
\def\Htwentythree{\lbrack \Height , C\rbrack}
\def\Htwentyfour{\lbrack \Hnine , C\rbrack}
\def\Htwentyfive{\lbrack \Hten , C\rbrack}
\def\Htwentysix{\lbrack \Heleven , C\rbrack}
\def\Htwentyseven{\lbrack \Htwelve , C\rbrack}
\def\Htwentyeight{\lbrack \Hthirteen , C\rbrack}
\def\Htwentynine{\lbrack \Hfourteen , C\rbrack}
\def\Hthirty{\lbrack \Hfive , \Hfour\rbrack}
\def\Hthirtyone{\lbrack \Hsix , \Hfour\rbrack}
\def\Hthirtytwo{\lbrack \Hsix , \Hfive\rbrack}

\def\LHfour{\lbrack  BA \rbrack}
\def\LHfive{\lbrack CA \rbrack}
\def\LHsix{\lbrack CB \rbrack}

\def\LHseven{\lbrack BA^2 \rbrack}
\def\LHeight{\lbrack CA^2 \rbrack}
\def\LHnine{\lbrack BAB \rbrack}
\def\LHten{\lbrack CAB \rbrack}
\def\LHeleven{\lbrack CB^2 \rbrack}
\def\LHtwelve{\lbrack BAC \rbrack}
\def\LHthirteen{\lbrack CAC \rbrack}
\def\LHfourteen{\lbrack CBC \rbrack}

\def\LHfifteen{\lbrack BA^3 \rbrack}
\def\LHsixteen{\lbrack CA^3 \rbrack}
\def\LHseventeen{\lbrack BA^2B \rbrack}
\def\LHeighteen{\lbrack CA^2B \rbrack}
\def\LHnineteen{\lbrack BAB^2 \rbrack}
\def\LHtwenty{\lbrack CAB^2 \rbrack}
\def\LHtwentyone{\lbrack CB^3 \rbrack}
\def\LHtwentytwo{\lbrack BA^2C \rbrack}
\def\LHtwentythree{\lbrack CA^2C \rbrack}
\def\LHtwentyfour{\lbrack BABC \rbrack}
\def\LHtwentyfive{\lbrack CABC \rbrack}
\def\LHtwentysix{\lbrack CB^2C \rbrack}
\def\LHtwentyseven{\lbrack BAC^2 \rbrack}
\def\LHtwentyeight{\lbrack CAC^2 \rbrack}
\def\LHtwentynine{\lbrack CBC^2 \rbrack}
\def\LHthirty{\lbrack \LHfive , \LHfour\rbrack}
\def\LHthirtyone{\lbrack \LHsix , \LHfour\rbrack}
\def\LHthirtytwo{\lbrack \LHsix , \LHfive\rbrack}

\def\LHthirtyfive{\lbrack BA^3B \rbrack}
\def\LHthirtysix{\lbrack CA^3B \rbrack}
\def\LHthirtyseven{\lbrack BA^2B^2 \rbrack}
\def\LHthirtyeight{\lbrack CA^2B^2 \rbrack}

\def\LHforty{\lbrack CAB^3 \rbrack}

\def\LHfortytwo{\lbrack BA^3C \rbrack}

\def\LHfortyfour{\lbrack BA^2BC \rbrack}
\def\LHfortyfive{\lbrack CA^2BC \rbrack}
\def\LHfortysix{\lbrack BAB^2C \rbrack}
\def\LHfortyseven{\lbrack CAB^2C \rbrack}

\def\LHfortynine{\lbrack BA^2C^2 \rbrack}

\def\LHfiftyone{\lbrack BABC^2 \rbrack}
\def\LHfiftytwo{\lbrack CABC^2 \rbrack}

\def\LHfiftyfour{\lbrack BAC^3 \rbrack}

\def\LHfiftyseven{\lbrack \LHseven , \LHfour\rbrack}
\def\LHfiftyeight{\lbrack \LHeight , \LHfour\rbrack}
\def\LHfiftynine{\lbrack \LHnine , \LHfour\rbrack}
\def\LHsixty{\lbrack \LHten , \LHfour\rbrack}
\def\LHsixtyone{\lbrack \LHeleven , \LHfour\rbrack}
\def\LHsixtytwo{\lbrack \LHtwelve , \LHfour\rbrack}
\def\LHsixtythree{\lbrack \LHthirteen , \LHfour\rbrack}
\def\LHsixtyfour{\lbrack \LHfourteen , \LHfour\rbrack}

\def\LHsixtyfive{\lbrack \LHseven , \LHfive\rbrack}
\def\LHsixtysix{\lbrack \LHeight , \LHfive\rbrack}
\def\LHsixtyseven{\lbrack \LHnine , \LHfive\rbrack}
\def\LHsixtyeight{\lbrack \LHten , \LHfive\rbrack}
\def\LHsixtynine{\lbrack \LHeleven , \LHfive\rbrack}
\def\LHseventy{\lbrack \LHtwelve , \LHfive\rbrack}
\def\LHseventyone{\lbrack \LHthirteen , \LHfive\rbrack}
\def\LHseventytwo{\lbrack \LHfourteen , \LHfive\rbrack}

\def\LHseventythree{\lbrack \LHseven , \LHsix\rbrack}
\def\LHseventyfour{\lbrack \LHeight , \LHsix\rbrack}
\def\LHseventyfive{\lbrack \LHnine , \LHsix\rbrack}
\def\LHseventysix{\lbrack \LHten , \LHsix\rbrack}
\def\LHseventyseven{\lbrack \LHeleven , \LHsix\rbrack}
\def\LHseventyeight{\lbrack \LHtwelve , \LHsix\rbrack}
\def\LHseventynine{\lbrack \LHthirteen , \LHsix\rbrack}
\def\LHeighty{\lbrack \LHfourteen , \LHsix\rbrack}

\def\filtr{\white}
\def\filtrALL{\{\filtr_n\}_{n\in\N}}
\def\qp{q+q^{-1}}
\def\qm{q-q^{-1}}
\def\qplus{(q+q^{-1})}
\def\qminus{(q-q^{-1})}
\def\lowgreeks{\alpha^r\beta^s\gamma^t}
\def\indset{\mathcal{I}}
\def\indsetnext{\mathcal{J}}
\def\indsetnextnext{\mathcal{K}}
\def\indsetL4{\indset^*}
\def\indsetLfive{\indsetnext^*}
\def\indsetzero{\indsetnextnext_0}
\def\indsetA{\indsetnextnext_1}
\def\indsetB{\indsetnextnext_2}

\def\indsetD{\indsetnextnext_4}
\def\indsetE{\indsetnextnext_5}
\def\Basis5{\mathcal{B}}
\def\Reps5{\mathcal{S}}
\def\HallB{\H}
\def\defId{\I}
\def\filtseq{\{\white_n\}_{n\in\N}}
\def\derwhite{\lbrack \white,\white\rbrack}
\def\Fker{\white\derwhite\white}
\def\Abar{\bar{A}}
\def\Bbar{\bar{B}}
\def\Cbar{\bar{C}}
\def\Fcomm{\F\left[\Abar,\Bbar,\Cbar\right]}
\def\abar{\bar{\alpha}}
\def\bbar{\bar{\beta}}
\def\gbar{\bar{\gamma}}
\def\Ombar{\bar{\Omega}}
\def\lowgreeksbar{\abar^r\bbar^s\gbar^t}
\def\AW{AW}
\def\AWfull{\AW:=\AW_q(a,b,c)}
\def\freeassoc3{\F\left< A,B,C\right>}

\begin{document}

\title{{\bf A Lie algebra related to the universal\\Askey-Wilson algebra}}
\author{{\scshape Rafael Reno S. Cantuba} \\
Mathematics Department\\
De La Salle University Manila\\
{ 2401 Taft Ave., Manila, 1004 Philippines}\\
{\it rafael{\_}cantuba@dlsu.edu.ph}}

\date{}

\maketitle

\begin{abstract}

Let $\F$ denote an algebraically closed field. Denote the three-element set by $\genset=\{A,B,C\}$, and let $\F\monoid$  denote the free unital associative $\F$-algebra on $\genset$. Fix a nonzero $q\in\F$ such that $q^4\neq 1$. The universal Askey-Wilson algebra $\white$ is the quotient space $\F\monoid/\defId$, where $\defId$ is the two-sided ideal of $\F\monoid$ generated by the nine elements $UV-VU$, where $U$ is one of $A,B,C$, and $V$ is one of
\begin{equation} \qplus A+\frac{qBC-q^{-1}CB}{\qm},\nonumber
\end{equation}
\begin{equation} \qplus B+\frac{qCA-q^{-1}AC}{\qm},\nonumber
\end{equation}
\begin{equation} \qplus C+\frac{qAB-q^{-1}BA}{\qm}.\nonumber
\end{equation}
Turn $\F\monoid$ into a Lie algebra with Lie bracket $\lbrack X,Y\rbrack = XY-YX$ for all $X,Y\in\F\monoid$. Let $\freeLie$ denote the Lie subalgebra of $\F\monoid$ generated by $\genset$, which is also the free Lie algebra on $\genset$. Let $\LieABC$ denote the Lie subalgebra of $\white$ generated by $A,B,C$. Since the given set of defining relations of $\white$ are not in $\freeLie$, it is natural to conjecture that $\LieABC$ is freely generated by $A,B,C$. We give an answer in the negative by showing that the kernel of the canonical map $\F\monoid\rightarrow\white$ has a nonzero intersection with $\freeLie$. Denote the span of all Hall basis elements of $\freeLie$ of length $n$ by $\freeLie_n$, and denote the image of $\sum_{i=1}^n\freeLie_i$ under the canonical map $\freeLie\rightarrow\LieABC$ by $\LieABC_n$. We show that the simplest nontrivial Lie algebra relations on $\LieABC$ occur in $\LieABC_5$. We exhibit a basis for $\LieABC_4$, and we also exhibit a basis for $\LieABC_5$ if $q$ is not a sixth root of unity. \\


\end{abstract}

\section{Introduction}

Let $\F$ be an algebraically closed field and fix a nonzero $q\in\F$ such that $q^4\neq 1$. Given $a,b,c\in\F$, the \emph{Askey-Wilson algebra} with parameters $a,b,c$ is the unital associative $\F$-algebra $\AWfull$ defined as having generators $A,B,C$ and relations
\begin{eqnarray} 
A+\frac{qBC-q^{-1}CB}{q^2-q^{-2}}=\frac{a}{\qp},\nonumber\\
B+\frac{qCA-q^{-1}AC}{q^2-q^{-2}}=\frac{b}{\qp},\nonumber\\
C+\frac{qAB-q^{-1}BA}{q^2-q^{-2}}=\frac{c}{\qp}.\nonumber
\end{eqnarray}
The algebra $\AW$ was introduced in \cite{Zhed} in order to desrcibe the Askey-Wilson polynomials \cite{Ask}. A wide range of applications of the Askey-Wilson algebra is discussed in \cite[Section~1]{UAW}. These applications include integrable systems, quantum mechanics, the theory of quadratic algebras, Leonard pairs and Leonard triples, and quantum groups. A central extension of the Askey-Wilson algebra $\AW$ is introduced in \cite{UAW}, which is called the \emph{universal Askey-Wilson algebra}.

\begin{definition}[\protect{\cite[Definition~1.2]{UAW}}]\label{UAWdef}
The \emph{universal Askey-Wilson algebra} is the unital associative $\F$-algebra, which we denote by $\white$, defined as having  generators $A,B,C$, and relations which assert that the following are central in $\white$:
\begin{eqnarray} 
A+\frac{qBC-q^{-1}CB}{q^2-q^{-2}},\label{central1}\\
B+\frac{qCA-q^{-1}AC}{q^2-q^{-2}},\label{central2}\\
C+\frac{qAB-q^{-1}BA}{q^2-q^{-2}},\label{central3}
\end{eqnarray}
where $q$ is a nonzero scalar that is not a fourth root of unity.
\end{definition}

Our main object of study is the Lie subalgebra $\LieABC$ of $\white$ generated by $A,B,C$. We show that a set of defining relations for $\white$ cannot be expressed in terms of Lie algebra operations only, and yet this does not imply that $\LieABC$ is freely generated by $A,B,C$. Denote the free unital associative $\F$-algebra on the three-element set $\genset=\{A,B,C\}$ by $\F\monoid$, and the free Lie algebra on $\genset$ by $\freeLie$. Recall that $\freeLie$ is the Lie subalgebra of $\F\monoid$ generated by $A,B,C$. We use the basis of $\freeLie$ which was introduced by Hall \cite{Hall}. Let us call the images of the Hall basis elements under the canonical map $\freeLie\rightarrow\LieABC$ as the standard Lie monomials of $\LieABC$. We show that the kernel of the canonical map $\F\monoid\rightarrow\white$ has a nonzero intersection with $\freeLie$. The generators $A,B,C$ are the standard Lie monomials of length $1$. The standard Lie monomials of lengths $\geq 1$ are constructed according to some rules, which we shall discuss in later sections. We show that the simplest Lie algebra relations on $\LieABC$ occur at length $5$, and we determine a maximal linearly independent set of standard Lie monomials of length at most 5.

\section{Preliminaries}

Let $\F$ be an algebraically closed field. Throughout, by an $\F$-algebra we mean a unital associative $\F$-algebra. Let $\mathfrak{A}$ be an $\F$-algebra. Recall that an anti-automorphism of $\mathfrak{A}$ is a bijective $\F$-linear map $\psi:\mathfrak{A}\rightarrow\mathfrak{A}$ such that $\psi(fg)=\psi(g)\psi(f)$ for all $f,g\in\mathfrak{A}$. We turn $\mathfrak{A}$ into a Lie algebra with Lie bracket $\lbrack f,g\rbrack=fg-gf$ for $f,g\in\mathfrak{A}$. 

Let $\N=\{0,1,2,\ldots\}$ denote the set of natural numbers. Given a nonzero $n\in\N$, let $\genset$ denote an $n$-element set. We shall refer to any element of $\genset$ as a \emph{letter}. For $t\in\N$, by a \emph{word of length} $t$ on $\genset$ we mean a sequence of the form
\begin{equation}
X_1X_2\cdots X_t,\label{assocform}
\end{equation}
where $X_i\in\genset$ for $1\leq i\leq t$. Given a word $W$ on $\genset$, denote the length of $W$ by $|W|$. The word of length $0$ will be denoted by $1$. Let $\monoid$ denote the set of all words on $\genset$. Given words $X_1X_2\cdots X_s$ and $Y_1Y_2\cdots Y_t$ on $\genset$, their \emph{concatenation product} is $$X_1X_2\cdots X_sY_1Y_2\cdots Y_t.$$

We now recall the free $\F$-algebra $\F\monoid$. The $\F$-vector space $\F\monoid$ has basis $\monoid$. Multiplication in the $\F$-algebra $\F\monoid$ is the concatenation product. We endow $\F\monoid$ with a symmetric bilinear form $(\ ,\ )$ with respect to which $\monoid$ is an orthonormal basis. For any $f\in\F\monoid$ and any word $W$, the coefficient of $W$ in $f$ is $(f,W)$. 

Given $n\in\N$, the subspace of $\F\monoid$ spanned by all the words of length $n$ is the $n$\emph{-homogenous component} of $\F\monoid$. Observe that $\F\monoid$ is the direct sum of all the $n$-homogenous components for $n\in\N$. If $f$ is an element of the $m$-homogenous component and $g$ is an element of the $n$-homogenous component, then $fg$ is an element of the $(m+n)$-homogenous component. It follows that the set of all $n$-homogenous components of $\F\monoid$ for all $n\in\N$ is a grading of $\F\monoid$.

The following notation will be useful. Let $W=X_1X_2\cdots X_t$ denote a word on $\genset$. We define $W\rev$ to be the word  $X_tX_{t-1}\cdots X_1$ on $\genset$.
Let $\theta$ denote the $\F$-linear map
\begin{eqnarray}
\theta:  \F\monoid & \rightarrow & \F\monoid,\nonumber\\
 W & \mapsto & (-1)^{|W|}W\rev,\label{thetadef}
\end{eqnarray}
for any word $W$. By \cite[p.~19]{Reut}, the map $\theta$ is the unique anti-automorphism of the $\F$-algebra $\F\monoid$ that sends $X$ to $-X$ for any letter $X$. 

Let $\freeLie$ denote the Lie subalgebra of the Lie algebra $\F\monoid$ generated by $\genset$. Following \cite[Theorem~0.5]{Reut}, we call $\freeLie$ the \emph{free Lie algebra on} $\genset$.

\begin{proposition}[\protect{\cite[Lemma 1.7]{Reut}}]\label{ReutProp} For $f\in\freeLie$, we have $\theta(f)=-f$.
\end{proposition}

We now recall the notion of a \emph{Lie monomial} on $\genset$. The set of all {Lie monomials} on $\genset$ is the minimal subset of $\F\monoid$ that contains $\genset$ and is closed under the Lie bracket. Observe that $0$ is a Lie monomial. Let $U$ be a Lie monomial. Then $U$ is an element of some $n$-homogenous component of $\F\monoid$. We define the \emph{length} of the Lie monomial $U$ to be $n$. Observe that $0$ has length $n$ for any $n\in\N$. Any nonzero Lie monomial has a unique length. Observe that the set of all Lie monomials of length $1$ is $\genset$. We now consider an ordering of Lie monomials.

\begin{definition}[\protect{\cite[p.~581]{Bonf}}]\label{LieOrderDef} Fix an ordering $<$ on $\genset$. Suppose that the set of all Lie monomials of lengths $1,2,\ldots , t-1$ have been ordered such that $U<V$ if the length of $U$ is strictly less than that of $V$. If $U,V$ both have length $t$, and can be written as $U=\lbrack X_1,Y_1\rbrack,V=\lbrack X_2,Y_2\rbrack$, then we compare $U,V$ using the following rules: 
\begin{enumerate}
\item If $Y_1\neq Y_2$, then $U<V$ iff $Y_1<Y_2$.
\item If $Y_1= Y_2$, then $U<V$ iff $X_1<X_2$.
\end{enumerate}
\end{definition}

We now introduce a basis for $\freeLie$ consisting of Lie monomials. 

\begin{proposition}[\protect{\cite[Theorem 3.1]{Hall}}]\label{HallProp} Let $\HallB$ be the set of Lie monomials such that $\genset\subset\HallB$, and that for any $U,V\in\HallB$, the Lie monomial $\lbrack U,V\rbrack$ is also in $\HallB$ whenever the following conditions hold.
\begin{enumerate}
\item $U>V$.
\item If $U=\lbrack X,Y\rbrack$ for some Lie monomials $X,Y$, then $Y\leq V$.
\end{enumerate}
Then $\HallB$ is a basis for $\freeLie$, often referred to as the \emph{Hall basis} of $\freeLie$.
\end{proposition}

\begin{example}\label{L4Ex} Suppose $\genset=\{A,B,C\}$ and $A<B<C$. Then the elements of $\HallB$ of length at most $4$ are:
\begin{eqnarray}
A,B,C,\ \Hfour,\ \Hfive,\ \Hsix,\ \Hseven,\ \Height,\ \Hnine,\nonumber\\
\Hten,\ \Heleven,\ \Htwelve,\ \Hthirteen,\ \Hfourteen,\nonumber\\
\Hfifteen,\ \Hsixteen,\ \Hseventeen,\ \Heighteen,\nonumber\\
\Hnineteen,\ \Htwenty,\ \Htwentyone,\ \Htwentytwo,\nonumber\\
\Htwentythree,\ \Htwentyfour,\ \Htwentyfive,\ \Htwentysix,\nonumber\\
\Htwentyseven,\ \Htwentyeight,\ \Htwentynine,\nonumber\\
\Hthirty,\ \Hthirtyone,\ \Hthirtytwo.\label{Hall4}
\end{eqnarray}
Observe that the above Lie monomials are listed according to the ordering in Definition~\ref{LieOrderDef}.
\end{example}

Given a Lie algebra $\mathfrak{L}$ and $x,y\in\mathfrak{L}$, recall the adjoint linear map $$\ad x:\mathfrak{L}\rightarrow\mathfrak{L}$$ that sends $y\mapsto\lbrack x,y\rbrack$. Denote an arbitrary word on $\genset$ by $W=X_1X_2\cdots X_t$. The \emph{Lie bracketing from left to right} is the linear map $\F\monoid\rightarrow\freeLie$ that sends $1\mapsto 0$ and sends the word $W$ into a Lie monomial according to the following rules:
\begin{enumerate}
\item If $|W|=1$, then $W\mapsto W$.
\item Suppose that the images of all words of length $<|W|$ have been defined. Denote the image of $X_1X_2\cdots X_{t-1}$ by $V$. Then $$W\mapsto \left(-\ad X_t\right) \left(V\right) = \lbrack V,X_t\rbrack. $$
\end{enumerate}
That is, $X_1X_2\cdots X_t\mapsto \lbrack\lbrack\lbrack X_1,X_2\rbrack,\cdots\rbrack, X_t\rbrack$ for $t\geq 2$. A Lie monomial that is an image of some word under Lie bracketing from left to right is said to be \emph{left-normed}.

\begin{notation}\label{leftNot} Given a word $W$, we denote the image of $W$ under Lie bracketing from left to right by $\lbrack W\rbrack$.
\end{notation}

\begin{example} With reference to Example~\ref{L4Ex}, we rewrite \eqref{Hall4} using Notation~\ref{leftNot}.
\begin{eqnarray}
A,B,C,\ \LHfour,\ \LHfive,\ \LHsix,\ \LHseven,\ \LHeight,\ \LHnine,\nonumber\\
\LHten,\ \LHeleven,\ \LHtwelve,\ \LHthirteen,\ \LHfourteen,\nonumber\\
\LHfifteen,\ \LHsixteen,\ \LHseventeen,\ \LHeighteen,\nonumber\\
\LHnineteen,\ \LHtwenty,\ \LHtwentyone,\ \LHtwentytwo,\nonumber\\
\LHtwentythree,\ \LHtwentyfour,\ \LHtwentyfive,\ \LHtwentysix,\nonumber\\
\LHtwentyseven,\ \LHtwentyeight,\ \LHtwentynine,\nonumber\\
\LHthirty,\ \LHthirtyone,\ \LHthirtytwo.\label{Hall4left}
\end{eqnarray}
\end{example}

Throughout, by an \emph{ideal} of an $\F$-algebra $\mathfrak{A}$ we mean a two-sided ideal of $\mathfrak{A}$. By a \emph{Lie ideal} of a Lie algebra $\mathfrak{L}$ we mean an ideal of $\mathfrak{L}$ under the Lie algebra structure. We now recall the notion of algebras having generators and relations (i.e., having a presentation). Denote the elements of $\genset$ by $G_1,G_2,\ldots, G_n$.

Let $f_1,f_2,\ldots ,f_m\in\F\monoid$ and let $I$ be the ideal of $\F\monoid$ generated by $f_1,f_2,\ldots ,f_m$. We define $\F\monoid/I$ as the $\F$-algebra with generators $G_1,G_2,\ldots ,G_n$ and relations $f_1=0,f_2=0,\ldots ,f_m=0$. The Lie subalgebra of $\F\monoid/I$ generated by $\genset$ is $\freeLie/\left(I\cap\freeLie\right)$.

Let $g_1,g_2,\ldots ,g_m\in\freeLie$ and let $J$ be the Lie ideal of $\freeLie$ generated by $g_1,g_2,\ldots ,g_m$. We define $\freeLie/J$ as the Lie algebra with generators $G_1,G_2,\ldots ,G_n$ and relations $g_1=0,g_2=0,\ldots ,g_m=0$. 

Suppose $\mathfrak{L}$ is a Lie algebra (over $\F$) generated by $\genset$. Then there exists an ideal $\mathcal{K}$ of $\freeLie$ such that $\mathfrak{L}=\freeLie/\mathcal{K}$. Let $\phi:\freeLie\rightarrow\freeLie/\mathcal{K}$ be the canonical Lie algebra homomorphism. Then the following span $\mathfrak{L}$: 
\begin{equation}\label{stdeq}
\phi(U),\ \ \ \mbox{for }U\in\HallB. 
\end{equation}
We call \eqref{stdeq} the \emph{standard Lie monomials} of the Lie algebra  $\mathfrak{L}$. Observe that the list of the standard Lie monomials of $\mathfrak{L}$ is identical to the list of elements of $\HallB$. This is because the Lie algebra homomorphism $\phi$ fixes generators. We order the list of standard Lie monomials of $\mathfrak{L}$ in a manner analogous to that given in Definition~\ref{LieOrderDef}.

\section{The universal Askey-Wilson algebra}

Hereon, let $\F$ be an algebraically closed field, and fix a nonzero $q\in\F$ such that $q^4\neq 1$. We fix $\genset=\{A,B,C\}$. Let $\F\monoid$ be the free associative algebra on $\genset$. We use the ordering $A<B<C$ to construct the Hall basis $\HallB$ of the free Lie algebra $\freeLie$ on $\genset$. Define the following elements of the free algebra $\F\monoid$.
\begin{eqnarray}
\alpha &:=& \qplus A + \frac{qBC-q^{-1}CB}{\qm},\label{alphadef}\\
\beta &:=& \qplus B + \frac{qCA-q^{-1}AC}{\qm},\label{betadef}\\
\gamma &:=& \qplus C + \frac{qAB-q^{-1}BA}{\qm}.\label{gammadef}
\end{eqnarray}
We also define the following Lie products in $\F\monoid$.
\begin{eqnarray}
r_0:=\lbrack A,\alpha\rbrack,  & r_3:=\lbrack B,\alpha\rbrack, & r_6:=\lbrack C,\alpha\rbrack,\nonumber\\
r_1:=\lbrack B,\beta\rbrack,  & r_4:=\lbrack C,\beta\rbrack, & r_7:=\lbrack A,\beta\rbrack,\nonumber\\
r_2:=\lbrack C,\gamma\rbrack,  & r_5:=\lbrack A,\gamma\rbrack, & r_8:=\lbrack B,\gamma\rbrack.\nonumber
\end{eqnarray}
Define $\defId$ as the ideal of $\F\monoid$ generated by $r_0,r_1,\ldots,r_8$. 

With reference to Definition~\ref{UAWdef}, we express $\white$ as a quotient space of $\F\monoid$, and as a consequence make explicit the defining relations of $\white$.

\begin{proposition} $\white=\F\monoid/\defId$.
\end{proposition}
\begin{proof} Recall $\white$ has relations which assert that each of \eqref{central1},\eqref{central2},\eqref{central3} commutes with every element of $\white$. Equivalently, each of \eqref{central1},\eqref{central2},\eqref{central3} commutes with every generator $A,B,C$. Observe that each of $\alpha,\beta,\gamma$ is a scalar multiple of \eqref{central1},\eqref{central2},\eqref{central3}, respectively. Then it suffices to define $\white$ as having nine defining relations of the form $\lbrack X,\delta\rbrack$, where $X\in\{ A,B,C\}$ and $\delta\in\{ \alpha,\beta,\gamma\}$. By the definition of $\defId$, we get the desired result.\qed
\end{proof}
We denote the images of $\alpha,\beta,\gamma$ under the canonical map $\F\monoid\rightarrow\white$ by the same symbols. 

\begin{proposition} $r_0,r_1,\ldots,r_8\notin\freeLie$.
\end{proposition}
\begin{proof} Let $\theta$ denote the $\F$-linear map in \eqref{thetadef}. It is routine to show that in the free algebra $\F\monoid$, we have $\theta(r_i)+r_i\neq 0$ for $0\leq i\leq 8$. Use Proposition~\ref{ReutProp}.\qed
\end{proof}

By a \emph{a word in $\white$} we mean the image of an element of $\monoid$ under the canonical map $\F\monoid\rightarrow\white$. Observe that the list of all words in $\white$ is identical to the list of all the words on $\genset$ in the free algebra $\F\monoid$ since the canonical map $\F\monoid\rightarrow\white$ is an $\F$-algebra homomorphism that fixes generators. We also preserve the ordering of generators $A<B<C$ in $\white$. By a \emph{$\white$-word}, we mean all elements of $\white$ of the form 
\begin{equation}\label{DwordForm}
W\lowgreeks 
\end{equation}
where $W$ is a word in $\white$, and $r,s,t\in\N$.

We now recall some properties of $\white$ as studied in \cite{UAW}. Let $U=X_1X_2\cdots X_t$ be a $\white$-word, where $X_i$ is either a generator of $\white$ or one of $\alpha,\beta,\gamma$ for $1\leq i\leq t$. Without loss of generality, we assume $U$ is of the form \eqref{DwordForm} since $\alpha,\beta,\gamma$ are central in $\white$. By an \emph{inversion} for $W$ we mean an ordered pair $(j,k)\in\N^2$ such that $1\leq j<k\leq t$ and $X_j,X_k\in\{A,B,C\}$ such that $X_j>X_k$. Any $\white$-word with no inversions is said to be \emph{irreducible}. For instance, $CABA$ has $4$ inversions and $CB^2A$ has $5$, while the $\white$-words $A^2BC,AB^2C$ are irreducible. The shortest words for which inversions exist are $BA,CA,CB$ and using \eqref{alphadef} to \eqref{gammadef}, the following hold in both $\F\monoid$ and $\white$.
\begin{eqnarray}
BA &=& q^2AB+q\qplus\qminus C-q\qminus\gamma,\label{BArel}\\
CA &=& q^{-2}AC-q^{-1}\qplus\qminus B+q^{-1}\qminus\beta,\label{CArel}\\
CB &=& q^2BC+q\qplus\qminus A-q\qminus\alpha.\label{CBrel}
\end{eqnarray}
Consider the word $CABA$, one of the $4$ inversions in which is caused by the first two letters $C,A$. Substituting for $CA$ using \eqref{CArel}, the result is a linear combination of $ACBA,B^2A,BA\beta$, each having fewer inversions than $CABA$.  

\begin{remark}\label{reduceRem} Following \cite[p.~7]{UAW} and \cite[Theorem~1.2]{Berg}, for any $\white$-word $W$, there exists a finite number of steps of substituting for inversions using \eqref{BArel} to \eqref{CBrel} such that the final result is a unique linear combination of irreducible $\white$-words. It follows that a basis for $\white$ consists of the vectors
\begin{eqnarray}
A^iB^jC^k\lowgreeks, \ \ i,j,k,r,s,t,\in\N.\label{ABCbasis}
\end{eqnarray}
\end{remark}

Given subspaces $H,K$ of $\white$, define $HK:=\Span\{hk\ |\ h\in H,k\in K\}$. If $K$ is a subspace of $H$, we say that a subspace $K'$ of $H$ is a \emph{complement of $K$ in $H$} whenever 
$$H=K+K'.\dirsum $$
We now recall a filtration for $\white$ as given in \cite[Section~5]{UAW}. This filtration is a sequence $\filtseq$ of subspaces of $\white$ defined by
\begin{eqnarray}
\white_0 & := & \F 1,\nonumber\\
\white_1 & := &  \white_0+\Span\{A,B,C,\alpha,\beta,\gamma\},\nonumber\\
\white_n & := & \white_1\white_{n-1},\ \ \ n\geq 1,\nonumber
\end{eqnarray}
and has the following properties for all $i,j\in\N$.
\begin{eqnarray}
\white_i &\subseteq& \white_{i+1},\label{filtsubset}\\
\white &=& \bigcup_{n\in\N} \white_n,\nonumber\\
\white_i\white_j &=& \white_{i+j}.\label{filtprod}
\end{eqnarray}
Given $n\in\N$, a basis for $\white_n$ consists of the vectors
\begin{eqnarray}
A^iB^jC^k\lowgreeks, \ \ i,j,k,r,s,t,\in\N,\ \ i+j+k+r+s+t\leq n,\label{filtnbasis}
\end{eqnarray}
while the following vectors form a basis for a complement of $\white_n$ in $\white_{n+1}$
\begin{eqnarray}
A^iB^jC^k\lowgreeks, \ \ i,j,k,r,s,t,\in\N,\ \ i+j+k+r+s+t= n+1.\label{ABCcompbasis}
\end{eqnarray}
We denote the span of the vectors \eqref{ABCcompbasis} by $\white_n^c$.

By \cite[Lemma~6.1]{UAW}, the following elements of $\white$ coincide and are central.
\begin{eqnarray}
qABC + q^{2} A^2 + q^{-2} B^2 +q^{2}  C^2-qA\alpha-q^{-1}B\beta-qC\gamma,\\
qBCA + q^{2} A^2 + q^{2} B^2 +q^{-2}  C^2-qA\alpha-qB\beta-q^{-1}C\gamma,\\
qCAB + q^{-2} A^2 + q^{2} B^2 +q^{2}  C^2-q^{-1}A\alpha-qB\beta-qC\gamma,\\
q^{-1}CBA + q^{-2} A^2 + q^{2} B^2 +q^{-2}  C^2-q^{-1}A\alpha-qB\beta-q^{-1}C\gamma,\\
q^{-1}ACB + q^{-2} A^2 + q^{-2} B^2 +q^{2}  C^2-q^{-1}A\alpha-q^{-1}B\beta-qC\gamma\label{OmACB},\\
q^{-1}BAC + q^{2} A^2 + q^{-2} B^2 +q^{-2}  C^2-qA\alpha-q^{-1}B\beta-q^{-1}C\gamma.\label{OmBAC}
\end{eqnarray}
Denote this element by $\Omega$, which is called in \cite{UAW} as the \emph{Casimir element of $\white$}. As shown in \cite[Section~7]{UAW}, we have other bases for $\white,\white_n$ (for $n\in\N$) that involve $\Omega$. First, the following vectors form a basis for $\white$.
\begin{eqnarray}\label{OmBasis}
A^iB^jC^k\Omega^l\alpha^r\beta^s\gamma^t,\ \ \ i,j,k,l,r,s,t\in\N,\ \ ijk=0.
\end{eqnarray}
Given $n\in\N$, a basis for $\white_n$ consists of the vectors
\begin{eqnarray}
A^iB^jC^k\Omega^l\alpha^r\beta^s\gamma^t,  & i,j,k,l,r,s,t\in\N,\ ijk=0,\  i+j+k+3l+r+s+t\leq n,\label{filtnbasisOm}
\end{eqnarray}
while the following vectors form a basis for a complement of $\white_n$ in $\white_{n+1}$.
\begin{eqnarray}
A^iB^jC^k\Omega^l\alpha^r\beta^s\gamma^t,  & i,j,k,l,r,s,t\in\N,\ ijk=0,\  i+j+k+3l+r+s+t= n+1.\nonumber
\end{eqnarray}

Recall that $\white$ is a Lie algebra with Lie bracket $\lbrack X,Y\rbrack :=XY-YX$ for $X,Y\in\white$. Denote the derived algebra of $\white$ by $\derwhite$, and the ideal of $\white$ generated by $\derwhite$ by $\Fker$. It follows that the Lie subalgebra of $\white$ generated by $A,B,C$ is $L:=\freeLie/(\defId\cap\freeLie)$. Given nonzero $n\in\N$, denote the span of all Hall basis elements of $\freeLie$ of length $n$ by $\freeLie_n$. Denote the image of $\sum_{i=1}^n\freeLie_i$ under the canonical map $\freeLie\rightarrow\LieABC$ by $\LieABC_n$. It follows that all standard Lie monomials of $\LieABC$ of length at most $n$ span $\LieABC_n$. 

\begin{proposition}\label{inclProp} $\LieABC\subseteq\F A+\F B+\F C +\Fker$.
\end{proposition}
\begin{proof} By the definition of $L$, we have $L\subseteq \F A+\F B+\F C + \derwhite$. Since $\white$ has a multiplicative identity, we have $\derwhite\subseteq\Fker$. From these we get the desired set inclusion. \qed
\end{proof}

\begin{proposition}\label{LcenterProp} If $q$ is not a root of unity, then $L$ has zero center.
\end{proposition}
\begin{proof} Let $q$ be not a root of unity. Suppose that $Z(\LieABC)$ has a nonzero element $f$. Since the generators $A,B,C$ of $\white$ are also in $\LieABC$, we have $Z(\LieABC)\subseteq Z(\white)$. By \cite[Corollary~8.3]{UAW}, $Z(\white)$ is generated by $\alpha,\beta,\gamma,\Omega$.  Observe that there exists a filtration subspace $\white_n$ such that $f\in\white_n$, and that  $\white_n\cap Z(\white)$ has a basis consisting of the vectors
\begin{equation}
\Omega^l\lowgreeks,\ \ \  l,r,s,t\in\N,\ \ 3l+r+s+t\leq n.\label{finwhite}
\end{equation}
Since $f$ is nonzero, there exists a nonzero $c\in\F$ and a vector $\Omega^w\alpha^x\beta^y\gamma^z$ in \eqref{finwhite} such that 
\begin{equation}\label{preimf}
f-c\Omega^w\alpha^x\beta^y\gamma^z = g,
\end{equation}
where $g$ is a linear combination of vectors in \eqref{finwhite} other than $\Omega^w\alpha^x\beta^y\gamma^z$. Let $\Fcomm$ denote the $\F$-algebra of polynomials in three mutually commuting indeterminates $\Abar,\Bbar,\Cbar$, with coefficients from $\F$. As shown in \cite[p.~17]{UAW}, there exists a unique surjective $\F$-algebra homomorphism $\Psi:\white\rightarrow\Fcomm$ with kernel $\Fker$ that sends \begin{equation}\label{commhom}
A\mapsto\Abar,\ B\mapsto\Bbar,\ C\mapsto\Cbar.
\end{equation}
Under this homomorphism, denote the images of $\alpha,\beta,\gamma,\Omega$ by $\abar,\bbar,\gbar,\Ombar$, respectively. As shown in \cite[Lemma~11.3]{UAW}, we have
\begin{eqnarray}
\abar & = & \qplus\Abar+\Bbar\Cbar,\\
\bbar & = & \qplus\Bbar+\Abar\Cbar,\\
\gbar & = &  \qplus\Cbar+\Abar\Bbar,\\
\Ombar & = & \qplus\Abar\Bbar\Cbar-\Abar^2-\Bbar^2-\Cbar^2.
\end{eqnarray} 
 It is routine to show that the vectors 
\begin{equation}
\Abar,\Bbar,\Cbar,\Ombar^l\lowgreeksbar,\ \ \  l,r,s,t\in\N,\ \ 3l+r+s+t\leq n,\label{finFcomm}
\end{equation}
are linearly independent in $\Fcomm$.
Observe also that by Proposition~\ref{inclProp}, $$f\in \F A+\F B+\F C + \ker \Psi.$$ Applying $\Psi$ to both sides of \eqref{preimf}, we have
\begin{equation}\label{postimf}
c_1\Abar +c_2\Bbar + c_3\Cbar-c\Ombar^w\abar^x\bbar^y\gbar^z = \bar{g},
\end{equation}
where $\bar{g}$ is a linear combination of the vectors in \eqref{finFcomm} except $\Ombar^w\abar^x\bbar^y\gbar^z $. We get a contradiction from \eqref{postimf}. Therefore, $Z(\LieABC)=0$. \qed
\end{proof}

We end this section by discussing some properties of $\white$ related to the group $\PSL2$. We denote by $\PSL2$ the free product of the cyclic group of order two and the cyclic group of order three \cite{Alp}. Let $\rho,\sigma$ denote the generators of $\PSL2$ such that $\rho^3=1$ and $\sigma^2=1$. By \cite[Theorem 3.1]{UAW}, the group $\PSL2$ acts faithfully on $\white$ as a group of automorphisms in the following way:

\begin{center}\begin{tabular}{c|ccc|ccc}
$u$ &  $\quad A\quad$ & $\quad B\quad$ & $C$  & $\alpha$ & $\beta$ & $\gamma$ \\
\hline
$\rho(u)$ & $B$ & $C$ & $A$ & $\beta$ & $\gamma$ & $\alpha$\\
\hline
$\sigma(u)$ & $B$ & $A$ & $ C+(q-q^{-1})^{-1}[A,B]$ & $\beta$ & $\alpha$ & $\gamma$
\end{tabular}
\end{center}
By \cite[Theorem~6.4]{UAW}, $\Omega$ is fixed by $\rho,\sigma$. It is routine to show that given $n\in\N$, the filtration subspace $\white_n$ is invariant under $\rho$.

\begin{proposition} The Lie algebra $\LieABC$ is $\PSL2$-invariant.
\end{proposition}
\begin{proof} Let $\tau\in\{\rho,\sigma\}$. It suffices to argue in the following way. Show that the images of the generators $A,B,C$ under $\tau$ are in $\LieABC$, and show that if the images of $f,g\in\LieABC$ under $\tau$ are in $\LieABC$, then so is the image of $\lbrack f,g\rbrack$ under $\tau$. By the above table, we are done with the first step. For the second step, assume that the images of $f,g\in\LieABC$ under $\tau$ are in $\LieABC$. Since $\tau$ is an $\F$-algebra automorphism, we have $\tau\left(\lbrack f,g\rbrack\right)=\tau\left(fg-gf\right) = \lbrack \tau(f),\tau(g)\rbrack \in\LieABC.$\qed 
\end{proof}

\section{$\LieABC$ is not free}

In this section, all computations are done in the free algebra $\F\monoid$. Our goal is to show that $\defId\cap\freeLie$ contains a nonzero element.

\begin{proposition} In the free algebra $\F\monoid$, 
\begin{eqnarray}
\frac{\LHfour}{q\qminus} - \qplus C &=& AB-\gamma,\label{BAfree}\\
\frac{\LHfive}{q^{-1}\qminus} -\qplus B &=& AC-\beta.\label{CAfree}
\end{eqnarray}
\end{proposition}
\begin{proof} Use \eqref{BArel},\eqref{CArel} to get \eqref{BAfree},\eqref{CAfree}, respectively.\qed
\end{proof}

\begin{proposition} In the free algebra $\F\monoid$, 
\begin{eqnarray}
\frac{\LHseven}{q^2\qminus^2} - \frac{\qplus\LHfive}{q\qminus} &=& A^2B+\qplus AC-A\gamma+\frac{r_5}{q\qminus},\label{BAAfree}\\
\frac{\LHtwelve}{\qplus\qminus^2} &=& B^2-A^2+\frac{A\alpha-B\beta+r_1}{\qp}.\label{BACfree}
\end{eqnarray}
\end{proposition}
\begin{proof} Apply $-\ad A$ to both sides of \eqref{BAfree}. The linear combination in the right side of the resulting equation contains $ABA$ which can be further simplified using \eqref{BArel}. From this, we get \eqref{BAAfree}. We show \eqref{BACfree} holds. Apply $-\ad C$ to both sides of \eqref{BAfree}. The result involves $ACB$ in the right side, which can be further simplified using \eqref{CBrel}. From this, we get \eqref{BACfree}.\qed
\end{proof}

\begin{definition} We define the following elements of $\freeLie$.
\begin{eqnarray}
H_0 & := & \frac{\LHthirtyone-\LHtwentyfour}{\qplus^2\qminus^2}+\frac{\LHseven}{\qm}-2\LHfive\\
I_0 &:=& \lbrack H_0,\LHfour\rbrack\label{I0def}
\end{eqnarray}
\end{definition}

\begin{lemma} In the free algebra $\F\monoid$, 
\begin{eqnarray}
\frac{\LHfour\alpha}{\qplus^2} = H_0 +\frac{r_1B-Ar_3}{\qplus^2} +\frac{q\qminus r_5}{\qp}.\label{H4afree}
\end{eqnarray}
\end{lemma}
\begin{proof} Apply $-\ad B$ to both sides of \eqref{BACfree}. The left side involves $BA^2$ which can be written uniquely as a linear combination of $A^2B,AC,A\gamma,B,\beta$ by repeated use of the relations \eqref{BArel},\eqref{CArel}. We get
\begin{eqnarray}
\frac{\lbrack BACB\rbrack}{q^2\qplus^2\qminus^3} &=& A^2B +q^{-3}(q^4+1)AC-A\gamma\nonumber\\
& &-q^{-2}\qplus\qminus B +q^{-2}\qminus\beta\nonumber\\
& & +\frac{r_5}{q\qplus}-\frac{\LHfour\alpha+Ar_3-r_1B}{q^2\qplus^2\qminus}\label{BACBfree}
\end{eqnarray}
In \eqref{BACBfree}, eliminate $A^2B,A\gamma$ using \eqref{BAAfree}, and eliminate $AC,B,\beta$ using \eqref{CAfree}. In the resulting equation, express all Lie monomials in terms of Hall basis elements. From this, we get \eqref{H4afree}.\qed
\end{proof}

\begin{lemma} In the free algebra $\F\monoid$, 
\begin{eqnarray}
I_0 & = &  \LHfour\frac{\lbrack r_0, B\rbrack-\lbrack r_3 , A\rbrack}{\qplus^2} -\frac{q\qminus\lbrack r_5,\LHfour\rbrack}{\qp}+\frac{\lbrack Ar_3-r_1B ,\LHfour\rbrack}{\qplus^2}.\label{I0free}
\end{eqnarray}
\end{lemma}
\begin{proof} Apply $-\ad\LHfour$ to both sides of \eqref{H4afree}. The resulting left side is 
\begin{equation}
\frac{\lbrack\LHfour\alpha ,\LHfour\rbrack}{\qplus^2},
\end{equation}
which can be simplified into $$\LHfour\frac{\lbrack r_0, B\rbrack-\lbrack r_3 , A\rbrack}{\qplus^2}$$ using $r_0=\lbrack A,\alpha\rbrack , r_3=\lbrack B,\alpha\rbrack$. The result is \eqref{I0free}.\qed
\end{proof}

\begin{theorem} The Lie algebra $\LieABC=\freeLie/(\defId\cap\freeLie)$ is not freely generated by $A,B,C$.
\end{theorem}
\begin{proof} Observe that if we write $I_0$ in terms of Hall basis elements, we have
$$I_0=\frac{\lbrack\LHthirtyone ,\LHfour\rbrack-\lbrack \LHtwentyfour ,\LHfour\rbrack}{\qplus^2\qminus^2} +\frac{\LHfiftyseven}{\qm}-2\LHthirty,$$
which, by the linear independence of the Hall basis elements, implies that $I_0\neq 0$. But by \eqref{I0free}, we have $I_0\in\defId\cap\freeLie$. Therefore, $\freeLie/(\defId\cap\freeLie)\neq \freeLie$.\qed
\end{proof}

\section{Properties of some standard Lie monomials of $\LieABC$}

We discuss properties of some standard Lie monomials of $\LieABC$ in relation to the filtration $\filtseq$ of $\white$.

\begin{proposition}\label{twogensProp0} For any $i,j\in\N$, the following hold in $\white$.
\begin{eqnarray}
\lbrack BA^i\rbrack-q^i\qminus^iA^iB & \in & \filtr_i,\label{LBAi}\\
\lbrack CA^i\rbrack- (-1)^iq^{-i}\qminus^iA^iC & \in & \filtr_i,\label{LCAi}\\
\lbrack CB^j\rbrack-q^j\qminus^jB^jC & \in & \filtr_j.\label{LCBj}
\end{eqnarray}
\end{proposition}
\begin{proof} We show \eqref{LBAi} holds by induction on $i$. The case $i=0$ is trivial. Suppose that for some $i\in\N$, we have 
\begin{equation}\label{LBAiminus1}
\lbrack BA^{i-1}\rbrack-q^{i-1}\qminus^{i-1}A^{i-1}B \in  \filtr_{i-1}.
\end{equation}
Denote the element in \eqref{LBAiminus1} by $X$. By the properties of the filtration $\filtrALL$, we have $\lbrack X,A\rbrack\in\filtr_i$. Using \eqref{LBAi}, we further obtain $$\lbrack X,A\rbrack+q^i\qminus^iA^{i-1}\left(\qplus C-\gamma\right)=\lbrack BA^i\rbrack-q^i\qminus^iA^iB,$$
which proves \eqref{LBAi}. The relations \eqref{LCAi} and \eqref{LCBj} are proven similarly.\qed
\end{proof}

\begin{proposition}[\protect{\cite[Lemma 8.1]{UAW}}]\label{TerProp1} Let $i,j,k\in\N$. Then the following hold in $\white$.
\begin{eqnarray}
\lbrack A,A^iB^jC^k\rbrack\  \ \ - & \left(1-q^{2(j-k)}\right)A^{i+1}B^{j}C^{k} & \in \  \ \filtr_{i+j+k},\label{AABC}\\
\lbrack B,A^iB^jC^k\rbrack\  \ \ - & \left(q^{2i}-q^{2k}\right)A^{i}B^{j+1}C^{k} & \in\  \ \filtr_{i+j+k},\label{BABC}\\
\lbrack C,A^iB^jC^k\rbrack\  \ \ - & \left(q^{2(j-i)}-1\right)A^{i}B^{j}C^{k+1} & \in  \ \ \filtr_{i+j+k}.\label{CABC}
\end{eqnarray}
\end{proposition}

\begin{proposition}\label{twogensProp1} For nonzero $i,j,k\in\N$, the following hold in $\white$.
\begin{eqnarray}
\lbrack BA^iB^j\rbrack-(-1)^jq^i(q^{2i}-1)^j\qminus^iA^iB^{j+1} & \in & \filtr_{i+j},\label{BAiBj}\\
\lbrack CA^iC^k\rbrack- (-1)^{i}q^{-i(2k+1)}(q^{2i}-1)^k\qminus^iA^iC^{k+1} & \in & \filtr_{i+k},\label{CAiCk}\\
\lbrack CB^jC^k\rbrack-(-1)^kq^j(q^{2j}-1)^k\qminus^jB^jC^{k+1} & \in & \filtr_{j+k}.\label{CBjCk}
\end{eqnarray}
\end{proposition}
\begin{proof} To show \eqref{BAiBj}, use the relation \eqref{LBAi}, the relation \eqref{BABC} with $k$ set to zero, and induction on $j$. The relations \eqref{CAiCk} and \eqref{CBjCk} are proven similarly.\qed
\end{proof}

\begin{proposition}\label{H10H12prop} The complement $\filtr_1^c$ of $\filtr_1$ in $\filtr_2$ contains $\LHten$ and $\LHtwelve$.
\end{proposition}
\begin{proof} Use the canonical map $\F\monoid\rightarrow\white$ on \eqref{BACfree} in order to obtain
\begin{eqnarray}
\frac{\LHtwelve}{\qplus\qminus^2}  =  B^2 - A^2 +\frac{A\alpha-B\beta}{\qp}\in\filtr_1^c.\label{BACincomp}
\end{eqnarray}
Apply $-\rho^2$ to both sides of \eqref{BACincomp}. We get
\begin{eqnarray}
\frac{\LHten}{\qplus\qminus^2}  =  C^2 - A^2 +\frac{A\alpha-C\gamma}{\qp}\in\filtr_1^c.\nonumber\qed
\end{eqnarray}
\end{proof}

\begin{proposition}\label{length5Prop0} For nonzero $i,j,k\in\N$ with $i\geq 2$, the following hold in $\white$.
\begin{eqnarray}
\lbrack BAB^jC^k\rbrack-  (-1)^{j+k}q^{-j}(q^{2j}-1)^k\qminus^{j+1}B^jC^{k-1}\Omega & \in & \filtr_{j+k+1},\ \ \ \ \ \  \label{BABjCk}\\
\lbrack CA^iB^j\rbrack -  (-1)^{i+j}q^{1-i}(q^{2(i-1)}-1)^j\qminus^{i}A^{i-1}B^{j-1}\Omega & \in & \filtr_{i+j},\label{CAiBj}\\
\lbrack CA^iBC^k\rbrack -  (-1)^{i+1}q^{(1-i)(1+2k)}(q^{2(i-1)}-1)^{k+1}\qminus^{i}A^{i-1}C^{k}\Omega & \in & \filtr_{i+k+1}.\label{CAiBCk}
\end{eqnarray}
\end{proposition}
\begin{proof} We show \eqref{BABjCk} holds. We first consider the case $k=1$.  By setting $i=1$ in \eqref{BAiBj}, we get
\begin{equation} \label{BABj}
\lbrack BAB^j\rbrack - (-1)^jq^{j+1}\qminus^{j+1}AB^{j+1}\in\filtr_{j+1}.
\end{equation}
Apply $-\ad C$ to the element in \eqref{BABj}. We get
\begin{equation} \label{BABj1pt5}
\lbrack BAB^jC\rbrack - (-1)^jq^{j+1}\qminus^{j+1}\lbrack AB^{j+1},C\rbrack\in\filtr_{j+2}.
\end{equation}
In \eqref{CABC} set $i,j,k$ to $1,j+1,0$, respectively, and combine with \eqref{BABj1pt5}. We have
\begin{equation} \label{BABj2}
\lbrack BAB^jC\rbrack - (-1)^{j+1}q^{j+1}(q^{2j}-1)\qminus^{j+1}AB^{j+1}C\in\filtr_{j+2}.
\end{equation}
Using \eqref{BArel}, it is routine to show that 
\begin{equation} \label{BABj3}
AB^n-q^{-2n}B^nA\in\filtr_{n},
\end{equation}
for $n\in\N$.  Set $n=j+1$ in \eqref{BABj3} and multiply the element by $C$ from the right. We get
\begin{equation} \label{BABj4}
AB^{j+1}C-q^{-2(j+1)}B^{j+1}AC\in\filtr_{j+2}.
\end{equation}
From \eqref{BABj2} and \eqref{BABj4}, we get 
\begin{equation} \label{BABj5}
\lbrack BAB^jC\rbrack -(-1)^{j+1}q^{-(j+1)}(q^{2j}-1)\qminus^{j+1}B^{j+1}AC\in\filtr_{j+2}.
\end{equation}
Using the fact that $\Omega$ is equal to \eqref{OmBAC}, we have
\begin{equation} \label{BABj6}
B^{j+1}AC-qB^j\Omega\in\filtr_{j+2}.
\end{equation}
From \eqref{BABj5} and \eqref{BABj6}, we get
\begin{equation} \label{BABj7}
\lbrack BAB^jC\rbrack -(-1)^{j+1}q^{-j}(q^{2j}-1)\qminus^{j+1}B^{j}\Omega\in\filtr_{j+2},
\end{equation}
from which we see that \eqref{BABjCk} holds for $k=1$ and for nonzero $j\in\N$. Using \eqref{CABC},\eqref{BABj7} and induction on $k$, we find that \eqref{BABjCk} holds for nonzero $j,k\in\N$. We now show \eqref{CAiBj} holds. Since $i\geq 2$, we can rewrite \eqref{BABjCk} changing the exponents $j,k$ to $i-1,j$, respectively. 
\begin{equation}\label{BABj8}
\lbrack BAB^{i-1}C^j\rbrack -  (-1)^{i+j-1}q^{1-i}(q^{2(i-1)}-1)^j\qminus^{i}B^{i-1}C^{j-1}\Omega \in  \filtr_{i+j}
\end{equation}
Denote the element in \eqref{BABj8} by $X$. Since $\filtr_{i+j}$ is invariant under $\rho$, we have $-\rho^2(X)\in\filtr_{i+j}$, where  $-\rho^2(X)$ is the element in \eqref{CAiBj}. Thus, \eqref{CAiBj} holds for nonzero $i,j\in\N$. Finally, we show \eqref{CAiBCk} holds. Set $j=1$ in \eqref{CAiBj}.
\begin{equation}\label{BABj9}
\lbrack CA^iB\rbrack -  (-1)^{i+1}q^{1-i}(q^{2(i-1)}-1)\qminus^{i}A^{i-1}\Omega \in \filtr_{i+1}
\end{equation}
Since $\Omega$ is central, if we apply $-\ad C$ to the element in \eqref{BABj9}, we get
\begin{equation}\label{BABj9pt5}
\lbrack CA^iBC\rbrack -  (-1)^{i+1}q^{1-i}(q^{2(i-1)}-1)\qminus^{i}\lbrack A^{i-1},C\rbrack\Omega \in \filtr_{i+2}.
\end{equation}
From \eqref{CABC} we obtain
\begin{equation}\label{BABj10}
\lbrack A^{i-1},C\rbrack\Omega -  q^{-2(i-1)}(q^{2(i-1)}-1)A^{i-1}C\Omega \in \filtr_{i+2}.
\end{equation}
From \eqref{BABj9} and \eqref{BABj10},
\begin{equation}\label{BABj11}
\lbrack CA^iBC\rbrack -  (-1)^{i+1}q^{(1-i)\cdot 3}(q^{2(i-1)}-1)^2\qminus^{i} A^{i-1}C\Omega \in \filtr_{i+2},
\end{equation}
from which we see that \eqref{CAiBCk} holds for $k=1$. Using \eqref{CABC},\eqref{BABj11} and induction on $k$, we find that \eqref{CAiBCk} holds for nonzero $i,k\in\N$ with $i\geq 2$.\qed
\end{proof}

\begin{proposition}\label{genOmProp} The following hold in $\white$.
\begin{eqnarray}
\LHeighteen & + \ \ \ \qminus^3A\Omega & \in \ \ \filtr_3,\label{Om18}\\
\LHtwentyfour & -\ \ \  \qminus^3B\Omega & \in \ \ \filtr_3,\label{Om24}\\
\LHthirtytwo & -\ \ \  \qminus^3C\Omega & \in\ \ \filtr_3.\label{Om32}
\end{eqnarray}
\end{proposition}
\begin{proof} The relations \eqref{Om18},\eqref{Om24} follow from \eqref{CAiBj},\eqref{BABjCk}, respectively. We show \eqref{Om32} holds. Let $V:=-\LHtwentyseven+\LHtwentyfive$. By Proposition~\ref{H10H12prop}, we have $V\in\filtr_3$. Denote the element in \eqref{Om24} by $X$. Using the fact that $\filtr_3$ is invariant under $\rho$, we have$$\rho(X) = \lbrack CBCA\rbrack -\qminus^3C\Omega\in\filtr_3. $$ Using the Jacobi identity to express $\lbrack CBCA\rbrack$ in terms of standard Lie monomials, we further have $$\rho(X) = \LHthirtytwo - V -\qminus^3C\Omega\in\filtr_3, $$ and it follows that
$$ \LHthirtytwo  -\qminus^3C\Omega=V+\rho(X) \in\filtr_3.\qed $$
\end{proof} 

\section{The standard Lie monomials of $\LieABC$ of length at most $4$}

Recall that the span of the standard Lie monomials of $\LieABC$ of length at most $n$ is $\LieABC_n$. Our goal in this section is to show that the standard Lie monomials of $\LieABC$ of length at most $4$ are linearly independent, and hence form a basis for $\LieABC_4$.

\begin{proposition}\label{twogensProp2} For nonzero $j,k\in\N$, the following hold in $\white$.
\begin{eqnarray}
\lbrack BAB^j\rbrack-(-1)^jq^{(j+1)}\qminus^{j+1}AB^{j+1} & \in & \filtr_{j+1},\label{BAiBj11}\\
\lbrack CAC^k\rbrack+ q^{-(k+2)}\qminus^{k+1}AC^{k+1} & \in & \filtr_{k+1},\label{CAiCk11}\\
\lbrack CBC^k\rbrack-(-1)^kq^{(k+1)}\qminus^{k+1}BC^{k+1} & \in & \filtr_{k+1},\label{CBjCk11}\\
\lbrack BA^2B^j\rbrack-(-1)^jq^{2(j+1)}\qplus^{j}\qminus^{j+2}A^2B^{j+1} & \in & \filtr_{j+2},\label{BAiBj12}\\
\lbrack CA^2C^k\rbrack- q^{-2(k+1)}\qplus^{k}\qminus^{k+1}A^2C^{k+1} & \in & \filtr_{k+2},\label{CAiCk12}\\
\lbrack CB^2C^k\rbrack-(-1)^kq^{2(k+1)}\qplus^{k}\qminus^{k+2}B^2C^{k+1} & \in & \filtr_{k+2}.\label{CBjCk12}
\end{eqnarray}
\end{proposition}
\begin{proof} Set $i=1,2$ in \eqref{BAiBj} to get \eqref{BAiBj11},\eqref{BAiBj12}. Do similarly to \eqref{CAiCk} and \eqref{CBjCk} in order to show the other relations.\qed
\end{proof}

\begin{lemma}\label{ind1Lem} Fix a nonzero $n\in\N$. The following vectors are linearly independent in $\white$ for any $i,j,k\in\N$ such that $1\leq i,j,k\leq n$.
\begin{eqnarray}
1,A,B,C,  \label{stateSTART} \\
\LHten,\LHtwelve,  \\
\LHeighteen,\LHtwentyfour, \LHthirtytwo,  \\
\lbrack BA^i\rbrack, \lbrack BAB^j\rbrack, \lbrack BA^2B^j\rbrack, \\
\lbrack CA^i\rbrack,  \lbrack CAC^k\rbrack, \lbrack CA^2C^k\rbrack,\\
\lbrack CB^j\rbrack,\lbrack CBC^k\rbrack,  \lbrack CB^2C^k\rbrack.\label{stateEND}
\end{eqnarray}
\end{lemma}
\begin{proof} Fix $n\in\N$. It suffices to show that there exists an upper triangular transition matrix from the above vectors to a subset of the basis of $\white$ consisting of the vectors in \eqref{OmBasis}:
\begin{equation}\nonumber
A^iB^jC^k\Omega^l\alpha^r\beta^s\gamma^t,\ \ \ i,j,k,l,r,s,t\in\N,\ \ ijk=0.
\end{equation}
Let $i,j,k\in\N$ such that $1\leq i,j,k\leq n$. From Propositions~\ref{twogensProp0}, \ref{H10H12prop}, \ref{genOmProp}, \ref{twogensProp2} we have the following data:
\begin{eqnarray}
\LHten -  c_1C^2 - d_1A^2 - d_2 C \gamma- d_3 A\alpha & \in & \filtr_0,\label{dataSTART}\\
\LHtwelve - c_2 B^2 - d_4 A^2 - d_5 B\beta - d_6 A\alpha & \in & \filtr_0,\label{dataotherEND}\\
\LHeighteen  - c_3 A\Omega & \in & \filtr_3,\label{dataotherSTART}\\
\LHtwentyfour  -  c_4 B\Omega & \in & \filtr_3,\\
\LHthirtytwo  -  c_5 C\Omega & \in & \filtr_3,\\
\lbrack BA^i\rbrack-e_{i}A^iB & \in & \filtr_{i},\\ 
\lbrack CA^i\rbrack- f_{i}A^iC & \in & \filtr_{i},\\
\lbrack CB^j\rbrack-g_{j}B^jC  & \in & \filtr_{j},\\
\lbrack BAB^j\rbrack-e_{j}'AB^{j+1} & \in & \filtr_{j+1},\\ 
\lbrack CAC^k\rbrack- f_{k}'AC^{k+1} & \in & \filtr_{k+1},\\
\lbrack CBC^k\rbrack-g_{k}'BC^{k+1}  & \in & \filtr_{k+1},\\
\lbrack BA^2B^j\rbrack-e_{j}''A^2B^{j+1} & \in & \filtr_{j+2},\label{conditionEx}\\ 
\lbrack CA^2C^k\rbrack- f_{k}''A^2C^{k+1} & \in & \filtr_{k+2},\label{conditionEx2}\\
\lbrack CB^2C^k\rbrack-g_{k}''B^2C^{k+1}  & \in & \filtr_{k+2},\label{dataEND}
\end{eqnarray}
where the small letters (other than $i,j,k$) denote scalars. Each of \eqref{dataotherSTART} to \eqref{dataEND} is of the form $M-aV \in\Delta_m$, where $M$ is a Lie monomial, $a\in\F$, and $V$ is an element of the basis of $\Delta$ consisting of the vectors in \eqref{OmBasis}, and $V\notin\Delta_m$. Call $V$ the \emph{leading term} of $M$. For \eqref{dataSTART},\eqref{dataotherEND}, define the leading terms of $\LHten$, $\LHtwelve$ by $C^2,B^2$, respectively. Observe that no two distinct Lie monomials found in \eqref{dataSTART} to \eqref{dataEND} have the same leading terms. This yields a transition matrix from the vectors \eqref{stateSTART} to \eqref{stateEND} to some of the vectors in \eqref{OmBasis} such that all entries below the main diagonal are zero, and that the diagonal entries are $$ c_1,\ldots, c_5, e_i,f_i,g_j,e_j',f_k',g_k',e_j'',f_k'',g_k'' . $$  By Propositions~\ref{twogensProp0}, \ref{H10H12prop}, \ref{genOmProp} and \ref{twogensProp2}, all such scalars are nonzero. Hence, the transition matrix is upper triangular.\qed
\end{proof}

\begin{notation} Let $\indset_n$ denote the set consisting of all the linearly independent vectors in Lemma~\ref{ind1Lem}.
\end{notation}

\begin{lemma}\label{ind2Lem} Fix nonzero $m,n\in\N$. The vectors $X\lowgreeks$ are linearly independent in $\white$ for any $X\in\indset_n$ and any $r,s,t\in\N$ such that $r+s+t\leq m$.
\end{lemma}
\begin{proof} The proof is similar to that of Lemma~\ref{ind1Lem}, but with \eqref{dataSTART} to \eqref{dataEND} modified as follows. For each of \eqref{dataSTART} to \eqref{dataEND}, multiply the element by $\lowgreeks$ and add $r+s+t$ to the index of the filtration subspace. Based on these new data, an upper triangular transition matrix can be constructed.\qed
\end{proof}

\begin{notation} Let $\indset_n^m$ denote the set consisting of all the linearly independent vectors in Lemma~\ref{ind2Lem}. Observe that the vectors
\begin{equation}\label{replace0}
\LHsix\gamma,\ \LHfour\beta,\ \LHfive\alpha,\ \LHsix\beta,\ \LHfive\gamma,\ \LHfour\alpha,
\end{equation}
are in $\indset_3^1$. Let $\indsetL4$ denote the set obtained from $\indset_3^1$ by replacing the vectors in \eqref{replace0} by the vectors 
\begin{equation}\label{replace1}
\LHtwenty,\ \LHtwentytwo,\ \LHtwentyfive,\ \LHtwentyseven,\ \LHthirty,\ \LHthirtyone.
\end{equation}
\end{notation}

\begin{proposition} The following hold in $\white$.
\begin{eqnarray}
\frac{\LHfour\alpha}{(q+q^{-1})^2} &=& \frac{\LHthirtyone-\LHtwentyfour}{(q+q^{-1})^2(q-q^{-1})^2}+\frac{\LHseven}{q-q^{-1}}-2\LHfive,\label{H4a}\\
\frac{\LHsix\beta}{(q+q^{-1})^2} &=& \frac{\LHtwentyseven-\LHtwentyfive}{(q+q^{-1})^2(q-q^{-1})^2}+\frac{\LHeleven}{q-q^{-1}}+2\LHfour,\label{H6b}\\
\frac{\LHfive\gamma}{(q+q^{-1})^2} &=& -\frac{\LHthirty+\LHeighteen}{(q+q^{-1})^2(q-q^{-1})^2}+\frac{\LHthirteen}{q-q^{-1}}-2\LHsix,\label{H5g}\\
\frac{\LHfour\beta}{(q+q^{-1})^2} &=& \frac{\LHthirty-\LHtwentytwo}{(q+q^{-1})^2(q-q^{-1})^2}-\frac{\LHnine}{q-q^{-1}}+2\LHsix,\label{H4b}\\
\frac{\LHsix\gamma}{(q+q^{-1})^2} &=& \frac{-\LHthirtyone+\LHtwentyfour-\LHtwenty}{(q+q^{-1})^2(q-q^{-1})^2}-\frac{\LHfourteen}{q-q^{-1}}-2\LHfive,\label{H6g}\\
\frac{\LHfive\alpha}{(q+q^{-1})^2} &=& \frac{-\LHtwentyfive}{(q+q^{-1})^2(q-q^{-1})^2}-\frac{\LHeight}{q-q^{-1}}-2\LHfour.\label{H5a}
\end{eqnarray}
\end{proposition}
\begin{proof} Apply the canonical map $\F\monoid\rightarrow\white$ to both sides of \eqref{H4afree} to  get \eqref{H4a}. Apply $\rho,\rho^2$ to  both sides of \eqref{H4a} to get \eqref{H6b},\eqref{H5g}, respectively. To get \eqref{H4b}, apply $\sigma$ to both sides of \eqref{H4a}. Apply $\rho,\rho^2$ to  both sides of \eqref{H4b} to get \eqref{H6g},\eqref{H5a}, respectively.\qed
\end{proof}

\begin{lemma}\label{L4Lem} The vectors in $\indsetL4$ are linearly independent in $\white$.
\end{lemma}
\begin{proof} We use \eqref{H4a} to \eqref{H5a} to construct a transition matrix from the elements of $\indset_3^1$ to the elements of $\indsetL4$. Denote such transition matrix by $T$. Order the rows of $T$ such that the last $17$ correspond to 
\begin{eqnarray}
& \LHfour,\  \LHfive,\ \LHsix,\ \LHseven,\ \LHeight,\ \LHnine,\label{aux1}\\
& \LHeleven ,\ \LHthirteen,\  \LHfourteen,\ \LHeighteen,\ \LHtwentyfour,\label{aux2}\\
& \LHtwenty,\ \LHtwentytwo,\ \LHtwentyfive,\  \LHtwentyseven,\ \LHthirty,\ \LHthirtyone,\label{replacements}
\end{eqnarray}
in that order, while order the columns of $T$ such that the last $17$ correspond to the vectors in \eqref{aux1},\eqref{aux2} together with
\begin{eqnarray}
\LHsix\gamma,\ \LHfour\beta,\LHfive\alpha,\ \LHsix\beta,\ \LHfive\gamma,\ \LHfour\alpha.\label{tobereplaced}
\end{eqnarray}
Observe that all the vectors in $\indset_3^1$ to be replaced to form $\indsetL4$ are in \eqref{tobereplaced}, all the replacements are in \eqref{replacements}, and all the other vectors that appear in \eqref{H4a} to \eqref{H5a} (which we use to construct the transition matrix) appear in \eqref{aux1},\eqref{aux2}. Then $T$ is of the form
\[
T=
\left[
\begin{array}{c|c}
I & M \\
\hline
0 & T'
\end{array}
\right]
\]
where $I$ is an identity matrix, $M$ is some matrix with $7$ columns, and $T'$ is a $7\times 7$ matrix which has the following properties. All diagonal entries of $T'$ are nonzero. Denote the $ij$-entry of $T'$ by $T'_{ij}$. All entries of $T'$ below the main diagonal and all entries in  the first two rows are zero except the ones that appear below:
$$\frac{1}{\qminus^2}=-T_{11}'=T_{31}'=-T_{71}'=-T_{22}'=T_{62}' \neq 0.$$ By these observations about $T'$, we find that $T$ is invertible. This implies that the vectors in $\indsetL4$ are linearly independent.\qed
\end{proof}

\begin{theorem} The standard Lie monomials of $\LieABC$ of length at most $4$ form a basis for $\LieABC_4$.
\end{theorem}
\begin{proof} All such vectors are in $\indsetL4$. Use Lemma~\ref{L4Lem} and the fact that the vectors in the statement span $\LieABC_4$. \qed
\end{proof}

\section{The standard Lie monomials of $\LieABC$ of length at most $5$}

In this section, we show the Lie algebra relations that hold in $\LieABC_5$. We also exhibit a basis for $\LieABC$ if $q$ is not a sixth root of unity. 

\begin{lemma} The following hold in $\white$.
\begin{eqnarray}
\LHfive\beta+\LHfour\gamma  & = & \frac{-\LHsixtyfive+\LHfiftyeight}{\qminus^3}-\frac{\LHtwentythree+\LHseventeen}{\qminus^2} \nonumber\\
& & -\frac{\qplus^2\left(\LHtwelve-\LHten\right)}{\qm}\label{H5b}\\
\LHsix\alpha-\LHfour\gamma & = & \frac{-\LHseventyfive+\LHsixtyone}{\qminus^3}-\frac{\LHtwentysix-\LHseventeen}{\qminus^2}\nonumber\\
& & -\frac{\qplus^2\LHten}{\qm}\label{H6a}
\end{eqnarray}
\end{lemma}
\begin{proof} In view of Remark~\ref{reduceRem}, write each of the left and right sides of \eqref{H5b} as a linear combination of irreducible $\white$-words. This yields the same linear combination of the basis vectors \eqref{ABCbasis} of $\white$. Apply $\rho$ to both sides of \eqref{H5b} to get \eqref{H6a}.\qed
\end{proof}

\begin{theorem}\label{rel5Thm} The following relations hold in $\LieABC$.
\begin{eqnarray}
\frac{\LHfortyfour}{\qm} &=& -\frac{-(2q^2+1)(q^2+2)\left(\LHseventythree+\LHsixtyseven\right)}{2q^2\qplus^2\qminus}\nonumber\\
& & -\frac{(q^4+3q^2+1)\left(\LHsixtytwo-2\LHsixty\right)}{2q^2\qplus^2\qminus}\nonumber\\
& & - \LHtwentyseven + 2\LHtwentyfive - \LHnineteen + \LHfifteen,\label{H44REL}\\
\frac{\LHsixtynine}{2\qplus^2\qminus} & = & \frac{-(3q^4+5q^2+3)\LHseventyeight}{2q^2\qplus^2\qminus}\nonumber\\
& & + \frac{(q^4+3q^2+1)\LHseventysix}{2q^2\qplus^2\qminus}\nonumber\\
& & - \frac{(2q^4+3q^2+2)\LHsixtyfour}{2q^2\qplus^2\qminus}\nonumber\\
& & -\frac{\LHfiftyone-\LHfortyseven}{\qm}\nonumber\\
& & +\LHtwentynine+\LHtwentytwo-\LHtwentyone+\LHeighteen,\label{H69REL}\\
\frac{\LHseventyfour}{2\qplus^2\qminus} & = & \frac{-2\LHseventy+\LHsixtyeight}{2q^2\qplus^2\qminus}\nonumber\\
& & + \frac{(2q^2+1)(q^2+2)\LHsixtythree}{2q^2\qplus^2\qminus}\nonumber\\
& & +\frac{\LHfortyfive}{\qm}-2\LHthirtyone - \LHtwentyeight \nonumber\\
& & + 2\LHtwentyfour - \LHtwenty + \LHsixteen,\label{H74REL}\\
\LHseventynine & = & -\LHseventyfive + \LHseventytwo - \LHsixtyfive \nonumber\\
& & + \LHsixtyone + \LHfiftyeight.\label{H79REL}
\end{eqnarray}
\end{theorem}
\begin{proof} To show \eqref{H44REL}, we first show that the equation
\begin{eqnarray}
\frac{\LHfortyfour}{\qplus^2\qminus^2} & = & \frac{-\LHtwentyseven + 2\LHtwentyfive}{\qm}\nonumber\\
& & +\frac{(q^4+1)\left(\LHnineteen-\LHfifteen\right)}{q^2\qminus}\nonumber\\
& & -\frac{2(q^6-1)\left(\LHeleven-\LHeight\right)}{q^3\qminus}\nonumber\\
& & +\frac{\LHtwelve\gamma}{\qplus^2}+\LHnine\beta+\LHseven\alpha\label{H44seed}
\end{eqnarray}
holds in $\white$. In view of Remark~\ref{reduceRem}, write each of the left and right sides of \eqref{H44seed} as a linear combination of irreducible $\white$-words. This yields the same linear combination of the basis vectors \eqref{ABCbasis} of $\white$. Apply $-\ad A,-\ad B,-\ad B,-\ad C$ to both sides of \eqref{H4a},\eqref{H5g},\eqref{H4b},\eqref{H6g}, respectively. Write all Lie monomials in standard form. We get
\begin{eqnarray}
\LHseven\alpha & = & f_1, \label{H44seed2}\\
\LHeight\beta & = & f_2, \label{H44seed3}\\
\LHten\gamma & = & f_3, \label{H44seed3pt5}\\
\LHtwelve\gamma +k\LHten\gamma & = & f_4, \label{H44seed4}
\end{eqnarray}
for some $k\in\F$ and some $f_1,f_2,f_3,f_4\in\LieABC$. Eliminate $\LHseven\alpha,\LHeight\beta,\LHtwelve\gamma$ in \eqref{H44seed} using \eqref{H44seed2} to \eqref{H44seed4}. The result is \eqref{H44REL}. Apply $\rho,\rho^2$ to both sides of \eqref{H44REL} in order to obtain \eqref{H69REL},\eqref{H74REL}, respectively. We now show \eqref{H79REL} holds. Add \eqref{H5b} and \eqref{H6a}. We get 
\begin{eqnarray}\label{g1g2}
\LHsix\alpha+\LHfive\beta = g_1 + g_2,
\end{eqnarray}
where $g_1,g_2$ are the right sides of \eqref{H5b},\eqref{H6a}, respectively.
Apply $-\rho$ to both sides of \eqref{H6a}. We get
\begin{eqnarray}\label{g3}
\LHsix\alpha + \LHfive\beta = g_3,\label{H6aDep}
\end{eqnarray}
for some $g_3\in\LieABC$ such that $\LHseventynine$ appears with nonzero coefficient in $g_3$. Eliminate $\LHsix\alpha + \LHfive\beta$ in \eqref{g1g2},\eqref{g3}. We get \eqref{H79REL} as desired.\qed
\end{proof}

At this point we have shown that each of 
\begin{eqnarray}\label{linDepsL5}
\LHfortyfour, \LHsixtynine, \LHseventyfour,\LHseventynine,
\end{eqnarray}
is linearly dependent on standard Lie monomials of length at most $5$ that are not in \eqref{linDepsL5}. In what follows, we shall show that the standard Lie monomials of length at most $5$ except \eqref{linDepsL5} are linearly independent in $\white$.

\begin{proposition} The following hold in $\white$.
\begin{eqnarray}
\LHthirtysix - \qplus\qminus^4A^2\Omega & \in & \filtr_4,\label{L5H36}\\
\LHthirtyeight-q\qminus^4AB\Omega & \in & \filtr_4,\\
\LHfortyfive+q^{-1}\qminus^4AC\Omega & \in & \filtr_4,\\
\LHfortysix+\qplus\qminus^4B^2\Omega & \in & \filtr_4,\\
\LHfiftyone+q\qminus^4BC\Omega & \in & \filtr_4,\label{L5H51}\\
\LHseventytwo +\qplus\qminus^4C^2\Omega & \in & \filtr_4.\label{L5H72}
\end{eqnarray}
\end{proposition}
\begin{proof} Use Proposition~\ref{length5Prop0} to show \eqref{L5H36} to \eqref{L5H51} hold. To show \eqref{L5H72}, set $i=2$ in \eqref{LCAi}. We get 
\begin{equation}\label{CCOm1}
\LHeight - q^{-2}\qminus^2A^2C\in\filtr_2.
\end{equation}
Applying $-\rho^2$ to the element in \eqref{CCOm1} and using the fact that $\filtr_2$ is invariant under $\rho$, we have
\begin{equation}\label{CCOm2}
\LHfourteen + q^{-2}\qminus^2C^2B\in\filtr_2.
\end{equation}
Set $i=1$ in \eqref{LCAi}. We get
\begin{equation}\label{CCOm3}
\LHfive + q^{-1}\qminus AC\in\filtr_1.
\end{equation}
By taking the Lie bracket of the elements in \eqref{CCOm2} and \eqref{CCOm3}, we have
\begin{equation}\label{CCOm4}
\LHseventytwo - q^{-3}\qminus^3\left(C^2BAC-AC^3B\right)\in\filtr_3\subset\filtr_4.
\end{equation}
Using \eqref{CArel}, it is routine to show that
$$C^2A-q^{-4}AC^2\in\filtr_2, $$
from which we obtain 
\begin{equation}\label{CCOm5}
\qminus^3\left(qC^2ACB-q^{-3}AC^3B\right)\in\filtr_4.
\end{equation}
Combining \eqref{CCOm4} and \eqref{CCOm5}, we obtain
\begin{equation}\label{CCOm6}
\LHseventytwo -\qminus^3C^2\left(q^{-3}BAC-qACB\right)\in\filtr_4.
\end{equation}
Using the fact that \eqref{OmACB} and \eqref{OmBAC} are both equal to $\Omega$, we have
\begin{equation}\label{CCOm7}
\qplus\qminus\Omega+q^{-3}BAC-qACB\in\filtr_2.
\end{equation}
Finally, we get \eqref{L5H72} from \eqref{CCOm6} and \eqref{CCOm7}.\qed
\end{proof}

\begin{proposition}\label{twogensPropL5} For nonzero $j,k\in\N$, the following hold in $\white$.
\begin{eqnarray}
\lbrack BA^3B^j\rbrack-(-1)^jq^3(q^{6}-1)^j\qminus^3A^3B^{j+1} & \in & \filtr_{j+3},\label{BAiBjL5}\\
\lbrack CA^3C^k\rbrack+ q^{-3(2k+1)}(q^{6}-1)^k\qminus^3A^3C^{k+1} & \in & \filtr_{k+3},\label{CAiCkL5}\\
\lbrack CB^3C^k\rbrack-(-1)^kq^3(q^{6}-1)^k\qminus^3B^3C^{k+1} & \in & \filtr_{k+3}.\label{CBjCkL5}
\end{eqnarray}
\end{proposition}
\begin{proof} Use Proposition~\ref{twogensProp1}.\qed
\end{proof}

\begin{lemma}\label{ind1LemL5} Assume $q$ is not a sixth root of unity. Fix a nonzero $n\in\N$. The following vectors are linearly independent in $\white$ for any $i,j,k\in\N$ such that $1\leq i,j,k\leq n$.
\begin{eqnarray}
1,A,B,C,  \label{L5stateSTART} \\
\LHten,\LHtwelve,  \\
\LHeighteen,\LHtwentyfour, \LHthirtytwo,  \\
\LHthirtysix, \LHthirtyeight, \LHfortyfive, \\
\LHfortysix, \LHfiftyone,\LHseventytwo,\\
\lbrack BA^i\rbrack, \lbrack BAB^j\rbrack, \lbrack BA^2B^j\rbrack, \lbrack BA^3B^j\rbrack,\\
\lbrack CA^i\rbrack,  \lbrack CAC^k\rbrack, \lbrack CA^2C^k\rbrack, \lbrack CA^3C^k\rbrack,\\
\lbrack CB^j\rbrack,\lbrack CBC^k\rbrack,  \lbrack CB^2C^k\rbrack, \lbrack CB^3C^k\rbrack.\label{L5stateEND}
\end{eqnarray}
\end{lemma}
\begin{proof} The proof is similar to that of Lemma~\ref{ind1Lem}. In order to construct the desired upper triangular transition matrix, we combine the data from \eqref{dataSTART} to \eqref{dataEND} to that in \eqref{L5H36} to \eqref{L5H72}, and \eqref{BAiBjL5} to \eqref{CBjCkL5}. Recall that the transition matrix that can be constructed from the data is upper triangular if the scalar coefficients of the leading terms are nonzero. Those in \eqref{dataSTART} to \eqref{dataEND} are nonzero as shown in the proof of Lemma~\ref{ind1Lem}. The scalar coefficients in \eqref{L5H36} to \eqref{L5H72} are nonzero by the manner $q$ is defined. Finally, the scalar coefficients in \eqref{BAiBjL5} to \eqref{CBjCkL5} are nonzero since $q$ is further assumed to be not a sixth root of unity.\qed
\end{proof}

\begin{notation} Let $\indsetnext_n$ denote the set consisting of all the linearly independent vectors in Lemma~\ref{ind1LemL5}.
\end{notation}

\begin{lemma}\label{ind2LemL5} Assume $q$ is not a sixth root of unity. Fix nonzero $m,n\in\N$. The vectors $Y\lowgreeks$ are linearly independent in $\white$ for any $Y\in\indsetnext_n$ and any $r,s,t\in\N$ such that $r+s+t\leq m$.
\end{lemma}
\begin{proof} The proof is similar to that of Lemma~\ref{ind1Lem}. For each of the data used in the proof of Lemma~\ref{ind1LemL5}, which are \eqref{dataSTART} to \eqref{dataEND}, \eqref{L5H36} to \eqref{L5H72}, and \eqref{BAiBjL5} to \eqref{CBjCkL5}, multiply the element by $\lowgreeks$ and add $r+s+t$ to the index of the filtration subspace. Use these new data to construct a similar upper triangular transition matrix.\qed
\end{proof}

\begin{notation} Let $\indsetnext_n^m$ denote the set consisting of all the linearly independent vectors in Lemma~\ref{ind2LemL5}. Observe that the vectors
\begin{equation}\label{replace2}
\LHsix\gamma,\ \LHfour\beta,\ \LHfive\alpha,\ \LHsix\beta,\ \LHfive\gamma,\ \LHfour\alpha,
\end{equation}
are in $\indsetnext_4^1$. Let $\indsetLfive$ denote the set obtained from $\indsetnext_4^1$ by replacing the vectors in \eqref{replace2} by the following vectors 
\begin{equation}\label{replace3}
\LHtwenty,\ \LHtwentytwo,\ \LHtwentyfive,\ \LHtwentyseven,\ \ \LHthirty,\ \ \LHthirtyone.
\end{equation}
\end{notation}

\begin{lemma}\label{L5Lem1} Assume $q$ is not a sixth root of unity. The vectors in $\indsetLfive$ are linearly independent in $\white$.
\end{lemma}
\begin{proof} The proof is similar to that of Lemma~\ref{L4Lem}.\qed
\end{proof}

\begin{proposition}\label{decompJ40Prop} Assume $q$ is not a sixth root of unity. If $V$ is a subspace of $\Span\indsetLfive$ such that $$\Span\indsetLfive=V+\Span\indsetnext_4^0 ,\dirsum$$ then a basis for $V$ is $\indsetLfive\backslash\indsetnext_4^0$.
\end{proposition}
\begin{proof} This follows from the fact that $\indsetLfive,\indsetnext_4^0$ are both linearly independent sets and that $\indsetnext_4^0\subset\indsetLfive$.\qed
\end{proof}

\begin{lemma}\label{replaceLem0} The following hold in $\white$.
\begin{eqnarray}
\LHthirteen\alpha + \frac{\LHfiftytwo}{\qminus^2} & \in & \Span\indsetnext_4^0\label{H52},\\
\LHfourteen\beta-\frac{\LHfiftyfour}{\qminus^2}+\frac{\LHfiftytwo}{\qminus^2}& \in & \Span\indsetnext_4^0\label{H54},\\
\LHseven\gamma-\frac{(q^6-1)\LHfiftyseven}{q^3\qminus^3}& \in & \Span\indsetnext_4^0\label{H57},\\
\LHnine\gamma+\frac{\qplus^2\LHfiftynine}{\qminus^2}& \in & \Span\indsetnext_4^0\label{H59},\\
\LHeight\beta-\frac{(q^6-1)\LHsixtysix}{q^3\qminus^3}& \in & \Span\indsetnext_4^0\label{H66},\\
\LHthirteen\beta+\frac{\qplus^2\LHseventyone}{\qminus^2} & \in & \Span\indsetnext_4^0\label{H71},\\
\LHeleven\alpha-\frac{(q^6-1)\LHseventyseven}{q^3\qminus^3} & \in & \Span\indsetnext_4^0\label{H77},\\
\LHfourteen\alpha+\frac{\qplus^2\LHeighty}{\qminus^2} & \in & \Span\indsetnext_4^0\label{H80}.
\end{eqnarray}
\end{lemma}
\begin{proof} Apply $-\ad C$ to both sides of \eqref{H5a},\eqref{H6b} and write all Lie monomials in standard form in order to get \eqref{H52},\eqref{H54}, respectively. To get \eqref{H57},\eqref{H59}, we first show that the equation
\begin{eqnarray}
\frac{\LHseven\gamma}{2\qplus^2} & = & \frac{(q^6-1)\LHfiftyseven}{2q^3\qplus^2\qminus^3}-\frac{\LHthirtyfive}{2\qplus^2\qminus^2}\nonumber\\
& & -\frac{\LHthirty+\LHtwentytwo}{2\qminus}+\LHthirteen-\LHnine\label{H7g}
\end{eqnarray}
holds in $\white$. In view of Remark~\ref{reduceRem}, write each of the left and right sides of \eqref{H7g} as a linear combination of irreducible $\white$-words. This yields the same linear combination of the basis vectors \eqref{ABCbasis} of $\white$. This proves \eqref{H7g}, from which \eqref{H57} follows. To prove \eqref{H59}, apply $\sigma$ to both sides of \eqref{H7g}. We obtain
\begin{eqnarray}
\frac{\LHnine\gamma}{2\qplus^2} & = & \frac{-\LHfiftynine}{2\qminus^2}-\frac{\LHthirtyseven}{2\qplus^2\qminus^2}\nonumber\\
& & +\frac{\LHthirtyone+\LHtwentyfour}{2\qminus}-\LHfourteen-\LHseven,\label{H9g}
\end{eqnarray}
from which \eqref{H59} follows. Finally, to prove \eqref{H66} to \eqref{H80}, apply $\rho,\rho^2$ to both sides of \eqref{H7g},\eqref{H9g}.\qed
\end{proof}

\begin{notation} Observe that the vectors
\begin{eqnarray}
 \LHthirteen\alpha,\ \LHfourteen\beta,\ \LHseven\gamma,\ \LHnine\gamma, \label{replace4a}\\
 \LHeight\beta,\ \LHthirteen\beta, \ \LHeleven\alpha,\  \LHfourteen\alpha,\label{replace4}
\end{eqnarray}
are in $\indsetLfive$. Let $\indsetzero$ denote the set obtained from $\indsetLfive$ by replacing the vectors in \eqref{replace4a},\eqref{replace4} by the vectors 
\begin{eqnarray}
\LHfiftytwo,\ \LHfiftyfour,\ \LHfiftyseven,\ \LHfiftynine,\\
\LHsixtysix,\ \LHseventyone,\ \LHseventyseven,\ \LHeighty.\label{replace5}
\end{eqnarray}
\end{notation}

\begin{lemma} Assume $q$ is not a sixth root of unity. The vectors in $\indsetzero$ are linearly independent in $\white$.
\end{lemma}
\begin{proof} We use Lemma~\ref{replaceLem0} in order to obtain a transition matrix from the vectors in $\indsetLfive$ into the vectors in $\indsetzero$. Denote each of \eqref{H52} to \eqref{H80} by $$M_i+f_i\in\Span\indsetnext_4^0,$$ where $1\leq i\leq 8$,  the vector $M_i$ is an element of $\indsetLfive$ that is to be replaced in order to form $\indsetzero$, while the standard Lie monomials that appear in $f_i$ are the replacements. Using the usual ordering of standard Lie monomials, define $\overline{M}_i$ as the largest standard Lie monomial in $f_i$. Let $j,k\in\N$, with $1\leq j,k\leq 8$. Observe that if $M_j\neq M_k$ then $\overline{M}_j\neq \overline{M}_k$. Observe also that $M_i\in\indsetLfive\backslash\indsetnext_4^0$ for $1\leq i\leq 8$. By Proposition~\ref{decompJ40Prop}, the coefficient of $M_i$ in $f_i$ is $-1$ for all $i$.  By these observations, it follows that there exists an upper triangular transition matrix from the  vectors in $\indsetLfive$ into the vectors in $\indsetzero$ with nonzero diagonal entries. These diagonal entries are precisely the scalar coefficients of $\overline{M}_i$ for all $i$. These coefficients are all nonzero since $q$ is assumed to be not a sixth root of unity. Since the vectors in  $\indsetLfive$ are linearly independent, the existence of a transition matrix just described implies that the vectors in $\indsetzero$ are also linearly independent.\qed
\end{proof}

\begin{theorem} Assume $q$ is not a sixth root of unity. The standard Lie monomials of $\LieABC$ of length at most $5$ except the vectors from \eqref{linDepsL5} form a basis for $\LieABC_5$.
\end{theorem}
\begin{proof} Let $\indsetnextnext$ denote the set obtained from $\indsetzero$ by replacing the vectors
\begin{eqnarray}
\LHtwelve\gamma,\ \LHten\gamma,\ \LHnine\beta,\ \LHseven\alpha,\label{removeA}\\
\LHseven\beta,\ \LHeight\gamma,\ \LHfive\beta,\label{removeB}\\
\LHtwelve\beta,\ \LHthirteen\gamma,\ \LHeight\alpha,\ \LHten\beta,\label{removeC}\\
\LHeleven\gamma,\ \LHsix\alpha,\ \LHnine\alpha,\label{removeD}\\
\LHten\alpha,\ \LHeleven\beta,\ \LHtwelve\alpha,\ \LHfourteen\gamma,\label{removeE}
\end{eqnarray}
by the vectors
\begin{eqnarray}
\LHsixty,\LHsixtytwo,\LHsixtyseven,\LHseventythree,\label{replaceA}\\
\LHfortytwo,\LHfiftyeight,\LHsixtyfive,\label{replaceB}\\
\LHfortynine,\LHsixtythree,\LHsixtyeight,\LHseventy,\label{replaceC}\\
\LHforty,\LHsixtyone,\LHseventyfive,\label{replaceD}\\
\LHfortyseven,\LHsixtyfour,\LHseventysix,\LHseventyeight.\label{replaceE}
\end{eqnarray}
We claim that $\indsetnextnext$ is linearly independent. Observe that all vectors mentioned in the statement of the theorem are in $\indsetnextnext$. By Theorem~\ref{rel5Thm}, the vectors in \eqref{linDepsL5} are linearly dependent on these vectors. The result follows. We now prove our claim.
We construct five more sets in a manner similar to the construction of $\indsetLfive$ from $\indsetnext_4^1$ and to that of $\indsetzero$ from $\indsetLfive$. The goal is that at each step, we prove that the constructed set is linearly independent. Let $\indsetA$ denote the set obtained from $\indsetzero$ by replacing the vectors \eqref{removeA} in $\indsetzero$ by the vectors \eqref{replaceA}. Let $\indsetB$ denote the set obtained from $\indsetA$ in a similar manner until $\indsetE$, which is obtained from $\indsetD$ by replacing the vectors \eqref{removeE} in $\indsetD$ by the vectors \eqref{replaceE}. Observe that $\indsetE=\indsetnextnext$. We show that each of $\indsetA,\ldots,\indsetE$ is a  linearly independent set in $\white$. Apply $-\ad A,-\ad B,-\ad B,-\ad A$ to both sides of \eqref{H6g},\eqref{H5g},\eqref{H4b},\eqref{H5a}, respectively. Write all Lie monomials in standard form. We get
\begin{eqnarray}
\LHtwelve\gamma +k\LHten\gamma  & = & f_1, \label{setAseed1}\\
\LHten\gamma & = & f_2, \label{setAseed2}\\
\LHnine\beta & = & f_3, \label{setAseed3}\\
\LHseven\alpha & = & f_4, \label{setAseed4}
\end{eqnarray}
for some $k\in\F$ and some $f_1,f_2,f_3,f_4\in\LieABC$. Eliminate $\LHfortyfour$ in \eqref{setAseed1},\eqref{setAseed3},\eqref{setAseed4} using the relation \eqref{H44REL}. Solve the resulting system in order to obtain
\begin{eqnarray}
\LHtwelve\gamma  & = & g_1, \label{setAseed1new}\\
\LHten\gamma & = & g_2, \label{setAseed2new}\\
\LHnine\beta & = & g_3, \label{setAseed3new}\\
\LHseven\alpha & = & g_4, \label{setAseed4new}
\end{eqnarray}
where each of $g_1,\ldots,g_4\in\LieABC$, is a linear combination of 
$$\LHsixty,\LHsixtytwo,\LHsixtyseven,\LHseventythree,$$
together with the vectors in $\indsetzero$. We use the equations \eqref{setAseed1new} to \eqref{setAseed4new} to construct a transtion matrix $T_1$ from the vectors in $\indsetzero$ into those in $\indsetA$. Index the columns of $T_1$ such that the last $4$ correspond to \eqref{removeA}, while index the rows such that the last $4$ correspond to \eqref{replaceA}. We find that $T_1$ has the form
\[
T_1=
\left[
\begin{array}{c|c}
I_1 & U_1 \\
\hline
0 & L_1
\end{array}
\right]
\]
where $I_1$ is an identity matrix, $U_1$ is some matrix with $4$ columns, and $L_1$ is a $4\times 4$ matrix, with 
$$\det L_1 =-\frac{1}{2\qplus^2\qminus^8}\neq 0,$$ which implies that $\det T_1\neq 0$. Thus, $\indsetA$ is linearly independent. For $2\leq i\leq 5$, we can also construct a transtion matrix $T_i$ from the vectors in $\indsetnextnext_{i-1}$ into those in $\indsetnextnext_i$ in a similar manner. The equations that can be used to construct $T_i$ are also derived from \eqref{H4a} to \eqref{H5a} with the application of the appropriate map $\ad X$ where $X\in\{A,B,C\}$. Furthermore, we find that $T_i$ can be partitioned into four matrices similar to $T_1$, and that if we denote the bottom right partition as $L_i$, then
\begin{eqnarray}
\frac{1}{\qminus^7}\ \ \  =& -\det L_2 = \det L_4 & \neq\ \  0, \nonumber\\
\frac{1}{\qminus^8}\ \ \  =& -\det L_3 = -\det L_5 & \neq \ \ 0,\nonumber
\end{eqnarray}
which imply that $T_i$ is invertible for $2\leq i\leq 5$. Therefore, $\indsetE=\indsetnextnext$ is a linearly independent set in $\white$.\qed
\end{proof}

\section{Acknowledgements}

This work was done while the author is a graduate student in De La Salle University, Manila, Philippines, and in part while he was an Honorary Fellow (Sept. to Dec. 2015) in the University of Wisconsin - Madison. He extends his thanks to his mentors, Prof. Arlene Pascasio and Prof. Paul Terwilliger. This research is supported by the Science Education Institute of the Department of Science and Techonology (DOST-SEI), Republic of the Philippines.

\end{document}